# An efficient co-simulation and control approach to tackle complex multi-domain energetic systems: concepts and applications of the PEGASE platform

*Mathieu Vallée[a], Roland Bavière[b], Valérie Seguin[b], Valéry Vuillerme[a], Nicolas Lamaison[a], Michaël Descamps[a] and Antoine Aurousseau[a]*

[a] *CEA LITEN, Le Bourget du Lac, France, firstname.lastname@cea.fr;*
[b] *CEA LITEN, Grenoble, France, firstname.lastname @cea.fr;*

**Abstract:** In this paper, we present a novel research software, called PEGASE, suitable for the design, validation and deployment of advanced control strategies for complex multi-domain energy systems. PEGASE especially features a highly efficient co-simulation engine, together with integrated solutions for defining both rule-based control strategies and Model-Predictive Control (MPC).
The main principle behind the PEGASE platform is divide-and-conquer. Indeed, rather than trying to solve a problem as a monolithic entity, which can be highly complex for multi-domain large-scale systems, it is often more efficient to decompose it into several domains or sub-problems, and to simulate them in a decoupled way. To provide its co-simulation capabilities, we based PEGASE on two main components. The first one is a framework for integrating simulation models, which can be either compatible with the FMI standard or interfaced through an Application Programming Interface (API). The second one is a multi-threaded sequencer enabling several simulation sequences with different time steps.
To provide advanced control capabilities, we also equipped PEGASE with a framework for MPC combining a comprehensive management of predictions data and a modeler dedicated to the formulation of Mixed Integer Linear Programs. We implemented this framework in C++ providing low formulation and resolution times for typical applications. Connection to hardware is also available via standard industry protocols thereby allowing PEGASE to control real energy systems.
In this paper, we show how these basic functionalities, combined with dedicated modeling tools, enable setting up simulation and control applications suitable for tackling the complexity of various kinds of energy systems. To illustrate this, we present four application examples from our recent research work. These examples cover several domains, from concentrated solar thermal plants to optimal control of district heating networks. The variety of examples demonstrates the robustness and genericity of the approach.

## 1. Introduction

Energy systems become more complex to define and operate in the context of ambitious environmental targets associated to an increase of the share of Renewable Energy Sources (RES).

To address these challenges, it is necessary to develop and integrate significant storage capacities, in order to introduce flexibility in the energy system and cope with demand-supply mismatches. Furthermore, energy systems consisting of multiple energy carrier allow converting an available form of energy into a required form of energy, depending on energy prices and on the storage/transportation feature of an energy carrier (e.g. electricity is easy to transport but difficult to store, while heat is easy to store and difficult to transport). However, taking advantage of these storage and conversion capacities often require fine-grained control relying on predictions, typically in the form of Model Predictive Control (MPC).

Numerical tools for simulation and optimization provide a very efficient way to design, optimize, operate and maintain the new generation of complex energy systems. In particular, accurate models (or "digital twins") enable the definition, evaluation and comparison of the energy systems as well as their control strategies in a controlled virtual environment, therefore extending the scope of exploration [1].

In this paper, we present a novel research software, called PEGASE, suitable for the design, validation and deployment of advanced control strategies for complex multi-domain energy systems. In section 2, we first give an overview of the current state of the art and need for efficient co-simulation tools. In section 3, we introduce the overall architecture and generic components of the PEGASE platform. In section 4, we present four different types of application designed using the PEGASE platform. We summarize this work in section 5.

## 2. State of the art

In the context of 0D-1D simulation, two families of tools can be identified. The first family encompasses domain specific tools, whereby the numerical methods are dedicated to a specific case of study. Such methods are usually computationally efficient even for large energy systems, however they are not suited for variants of the initial case, in particular if some physical mechanisms need to be customized. Examples of such tools include PowerFlow for electrical grids, Termis Engineering for heating networks, Energy Plus for building simulations. Examples for power plant include two-phase flow simulations with homogeneous model (3 equation- model) and with two-fluid model (5 or 6 equation-model) for nuclear applications (e.g. CATHARE).

The second family encompasses generic simulation platform, suitable for multi-physics and multi-domain simulations. Such platforms are flexible in the type of energy system which can be modelled, but to the detriment of the scale and robustness. This means that it can be more challenging to simulate large energy systems than with domain specific tools. Examples of such platform based on Modelica language includes DYMOLA, Simulation X or Open MODELICA. Other example relevant to energy systems include Simcenter Amesim and Matlab-Simulink.

Co-simulation platforms can be used to leverage the advantages of the different types of tools [2]. In practice, a complex simulation problem is divided into subsystems, which are easy to solve by the right tool, either a domain specific tool or a generic simulation platform. The robust coupling of the subsystems requires an appropriate implementation in order to avoid numerical instabilities and excessive run times. Functional Mock-up Interface (FMI) [3] is a tool independent standard to support co-simulation of dynamic models using a combination of xml-files and compiled C-code. Among the co-simulation platforms currently available, BCVTB (Building Controls Virtual Test Bed – Berkeley) [4] is an open source platform dedicated to building energy and control systems, whereas DACCOSIM (Distributed Architecture for Controlled CO-SIMulation) [5] is an open source and multi-physics platform.

When operating a complex energy system, the influence of the control strategy is particularly significant if energy storage components as well as multiple sources of energy with variable prices or intermittent renewable energy sources are present. Operational optimization approaches relying on various mathematical techniques have been proposed over the years (see for instance [6,7]) in the context of District Heating Systems). In particular, Model-Predictive Control (MPC) is an attractive solution. MPC is generally composed of a load prediction module and an optimization procedure used to determine the optimal trajectory for control variables, i.e. the one minimizing an objective function while meeting different technical and operational constraints.

The Mixed Integer Linear Programming (MILP) framework is one of the preferred solutions for formulating optimisation problems in the energy domain [8]. The first advantage is the close relation of a MILP formulation to the physical equations, which facilitates the expression of realistic models. The second advantage is the reasonable computational time, especially when using efficient modern MILP solvers. The third advantage is the warranty of optimality, which especially ensures the perfect reproducibility of results when performing comparisons.

In the perspective of an application to energy system management, a co-simulation platform should exhibit the following two characteristics: (i) it needs to be accessible to non-expert users, especially regarding the implementation of the sub-systems coupling, and (ii) it should enable the two types of control: Rule Based Control (RBC) and Model-Predictive Control (MPC). The PEGASE platform has been developed to satisfy these two objectives.

## 3. The PEGASE co-simulation and control platform

In this section, we introduce the architecture and generic components of the PEGASE platform. Fig. 1 presents the overall architecture of a typical application based on the PEGASE platform. It is organised in four layers, with some generic components and some application-specific components. Each of the generic components is described in the following, while section 4 gives some examples of application-specific modules.

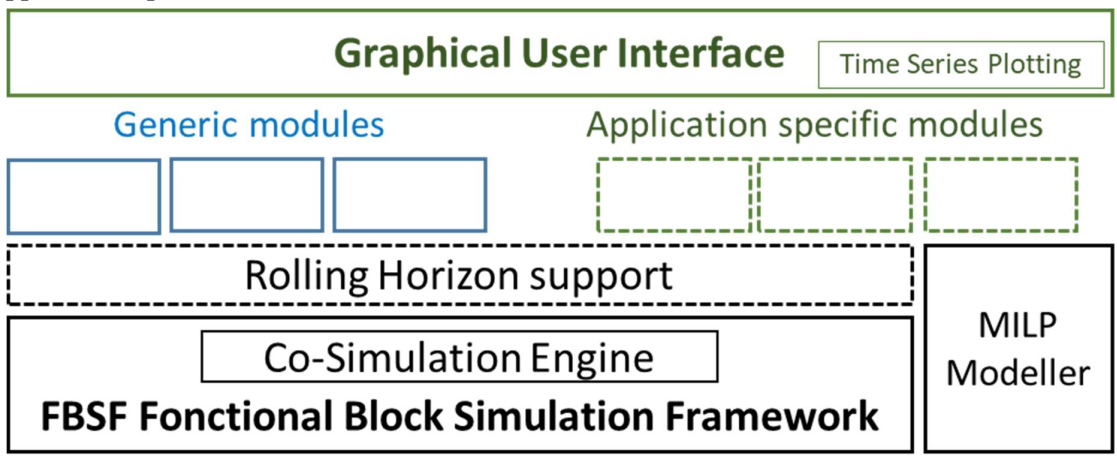

Figure 1: Architecture of a typical application based on the PEGASE platform.

### 3.1. The co-simulation engine

The base layer of the architecture is the commercial FBSF Functional Bloc Simulation Framework [9]. FBSF provides a very efficient co-simulation engine developed in C++ and fully compatible with the FMI 2.0 standard [3]. This co-simulation engine follows the domain decomposition approach. In this view, the global problem is decomposed into a collection of non-overlapping domains and each domain, or sub-model, is assigned to an ad hoc simulation tool. The coupling occurs by imposing that all computational variables are continuous at the domains' interfaces. In theory, an exact solution to the global simulation problem is only reachable using an implicit time marching scheme. It is also well known that an explicit time marching scheme can lead to numerical instability of the coupled calculation. However, for simplicity reasons and because these limitations can be circumvented in many practical situations, we programmed an explicit co-simulation engine in the PEGASE platform. This choice brings the two following remarks. First, a reasonable user of the PEGASE platform will conduct a time-step convergence study, i.e. a verification that the obtained results are independent of the coupling time-step. This guarantees that the potential inaccuracy of explicit coupling is not at play. Second, the explicit nature of the co-simulation engine requires a careful positioning of the boundaries of the computational domains. Yet, in practice, excluding a few well-known problematic situations such as breaking down a closed incompressible hydraulic loop into several sub-models is sufficient to limit the occurrence of coupling instability.

Each sub-model is either respecting the FMI standard or a dedicated code respecting a similar API. The supported methods allow for managing the life-cycle of each sub-model (e.g. construction, destruction ...), the time advance (e.g. pre-& post-stepping, time-step computation ...) and the data exchange through "mutators" and "accessors". Data are exchanged via a common "exchange zone", defined as the concatenation of the outputs of each sub-models. In the PEGASE platform, the execution of the sub-models is organized through a collection of sequences. The execution periodicity of the sequences all derive from a unique period and are constant over time. The multi-period character of the co-simulation engine is convenient to manage different coupling time steps or to trigger specific calculations required by the graphical user interface at low frequency.

## 3.2. Rolling Horizon and MILP Modeller support for Model-Predictive Control

The second layer of the architecture provides an optional support for programming Model-Predictive Control (MPC) strategies. MPC is one of the key approaches enabling the development of efficient Energy Management Systems (EMS). PEGASE integrates support for building MPC algorithms based on a rolling horizon principle. More specifically, it performs the following operations at regular time intervals:

1. Collect data available either from simulation models or from a real system, as well as weather and energy cost predictions. In particular, we consider systems with varying energy costs over the day, because of renewable energy usage or by considering multi energy sources with an optimal dispatch strategy.
2. Formulate an optimization problem, which generally consists in minimizing operational costs while respecting technical and non-technical constraints.
3. Solve the optimization problem over a given horizon (e.g. next 24h) in order to define the optimal trajectory of the control variables.
4. Apply the obtained set point for the next time interval, before performing the optimization again to adjust for prediction updates and real system behavior.

Such a control loop is enabled in PEGASE by manipulating so-called "time-vectors", i.e. multiple data points corresponding to the same data at various points in time. Indeed, the simulation and control modules exchange time-vectors in MPC mode, as opposed to single scalars in standard co-simulation. Additional support is provided for representing these "time-vectors" in time series graphs, and for managing time shifting between successive iterations.

The second aspect of the MPC support in PEGASE is the formulation and resolution of optimisation problems. In particular, many applications perform optimisation in the MILP framework, which has several advantages, as pointed out in section 2. The main drawback is the requirement of using only linear or quasi-linear models, which may require a simplification of non-linear models often used for simulating energy systems.

To ease the process of formulating MILP models in the context of energy systems, the PEGASE platform embeds MILP formulation capacities based on the Eigen linear algebra library [10]. The defined primitives allow formulating optimisation problems in a natural way, closely related to the mathematical formulation of the problem. Acquisition of data from the co-simulation as well as integration with numerous MILP solvers is also available. Using this functionality, problem formulation is performed at each iteration step in a few milliseconds on a standard PC, even for large problems (over 10k variables). It is important to notice that this feature especially allows formulating a different MILP problem at each iteration step, depending on the current state of the system under control. This is particularly useful for handling non-linearities in a flexible way.

### 3.3. Application modules

The third layer of the architecture is composed of the various modules used by a given application. Modules can be either external modules respecting the FMI standard, dedicated modules developed in C++ or so-called "logic" modules created using the graphical logic editor (see section 3.3.1**Erreur ! Source du renvoi introuvable.**). Several generic modules are available for common tasks, such as interacting with data sources (files, databases, industrial protocols), integrating external domain specific-codes (e.g. Cathare [11]), controlling specific aspects of the simulation (automatic stop, management of the MPC mode for simulation modules), or invoking external models in Python, Scilab or C++. For a given application, in addition to FMUs or user-defined "logic" modules, specific C++ modules can be developed to bring a greater level of flexibility. This is, however, not required for most common cases.

### 3.3.1. The graphical logic editor

The PEGASE software comes with a graphical editing and simulation framework to model causal and logic problems. This feature is particularly convenient to represent control-command and Rule-Based Control (RBC). The editor relies on model libraries containing generic components (e.g. PID regulator, logical switch …) which can be drag-and-dropped into a multi-page view of the global model.

A specificity of our causal simulation framework is that all model variables are accessible to visualization during run-time, just like in a debugger. This is very helpful during the development and verification phases. To illustrate this aspect, Fig. 2 shows a screenshot of a portion model during two consecutive time steps.

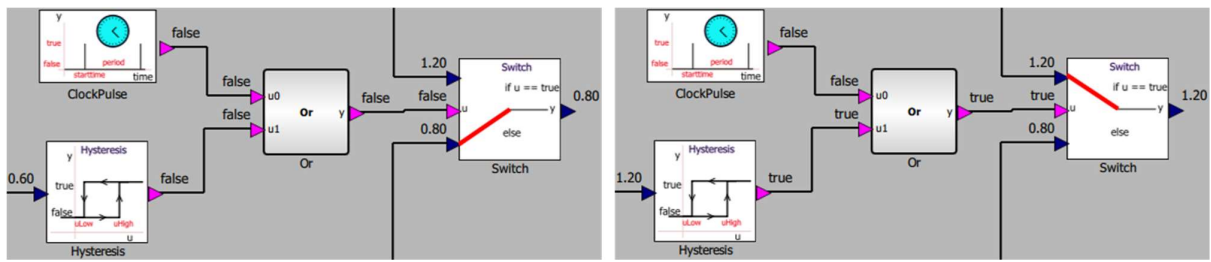

*Figure 2: Screenshot of a portion of logic model during run-time (left to right: 2 consecutive steps).*

The causal solver we programmed relies on a directed (dependency) graph, built during initialization, which links the inputs and outputs of all components of the model. The idea behind this formalization is that a directed acyclic graph, submitted to a topological sort, will directly provide a correct valid sequence of computational tasks. The question is then to transform the dependency graph into a directed acyclic graph. This is done by removing a feedback arc set, using a collection of one time-step delays.

### 3.4. Graphical user interface

The fourth layer is the graphical user interface, which consists in two parts. First, FBSF provides classical features for plotting time series, with a focus on online analysis. To this end, graphs are constantly updated and rescaled automatically while the calculations are running, and the view can be frozen to perform further analysis without stopping the underlying simulation or real system operation. Second, FBSF gives the possibility of programming elaborated graphical user interface in the Qt/QML framework [12]. When programming such user interfaces, all Qt/QML features can be used, and a dynamic connection to the underlying simulation data is available. This again gives great flexibility in the way simulation and control solutions can be presented in real-world application, as demonstrated for instance in large-scale district heating networks [13].

## 4. Sample PEGASE applications

In this section, we present several applications based on the PEGASE platform. These applications have been selected to demonstrate various capabilities. The first application (section 4.1) focuses on co-simulation combining two dedicated modelling tools. The second application (section 4.2) presents the generic methodology adopted for designing and validating the design of advanced Energy Management Systems (EMS). The third application (section 4.3) illustrates the deployment of such an EMS on a real-world demonstrator. Finally, the fourth application introduces the use of these methodologies for techno-economical assessment of an energy system (section 4.4).

### 4.1. Development and application of a concentrated solar power plant numerical simulator

The first application considers a numerical simulator for a Concentrated Solar Power (CSP) plant. Over the past years, CEA has been involved in the design, the realization and the operation of several CSP plants. Our research group has been particularly active on the so-called Direct Steam Generation variant (DSG), in which steam is directly generated in the solar field absorber tubes. The main challenge for the DSG-variant is to master the two-phase flow system under dynamic excitations, such as sun/clouds alternations [14]. In this context, the need for robust, reliable and validated simulation tool is obvious.

To meet this need, we developed a multi-domain dynamic simulator [15] of the system presented in Fig. 3. We decomposed the global simulation problem into: 1) a two-phase flow thermo-hydraulic problem we modelled using the domain specific 6-equations CATHARE code [11], 2) an optical and thermal problem for the mirror-absorber system programmed using the MODELICA language [16] and 3) several RBC models accounting for the low-level regulators as well as the safety functions of the CSP plant.

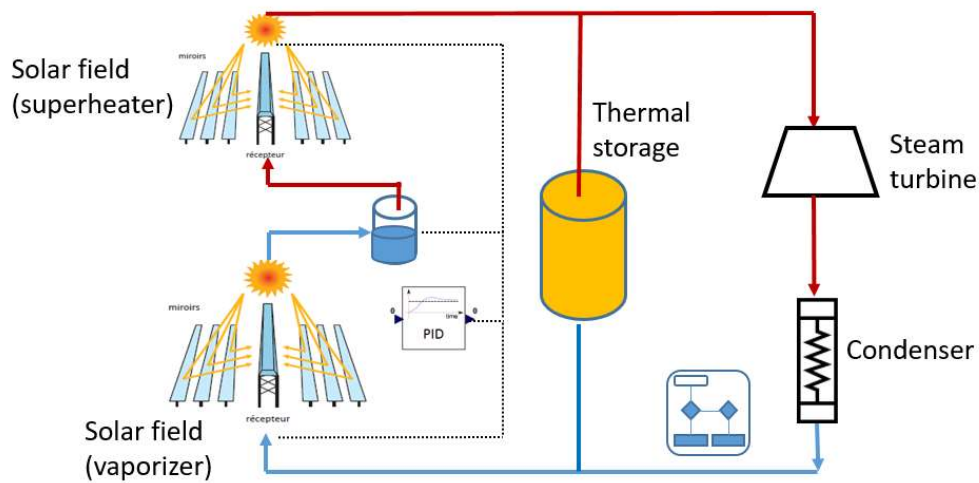

*Figure 3 : Scheme of the CSP prototype modelled using PEGASE.*

We positioned the boundary between the first two aforementioned computational domains at the outside surface of the absorber tubes. We dealt with the thermal coupling using a classical (temperature – heat-flux) explicit scheme. We computed the system within the PEGASE platform using a unique time-step of 2 seconds. This value is an upper limit both for the correct simulation of a sampled regulator and for the CATHARE code. These choices did not lead to any numerical instability.

We used the multi-domain simulator to assess various control strategies. For example, Fig. 4 shows the evolution of a relevant plant variable during a typical 15 minutes cloud shadow, for two different regulators implemented in the graphical logic editor.

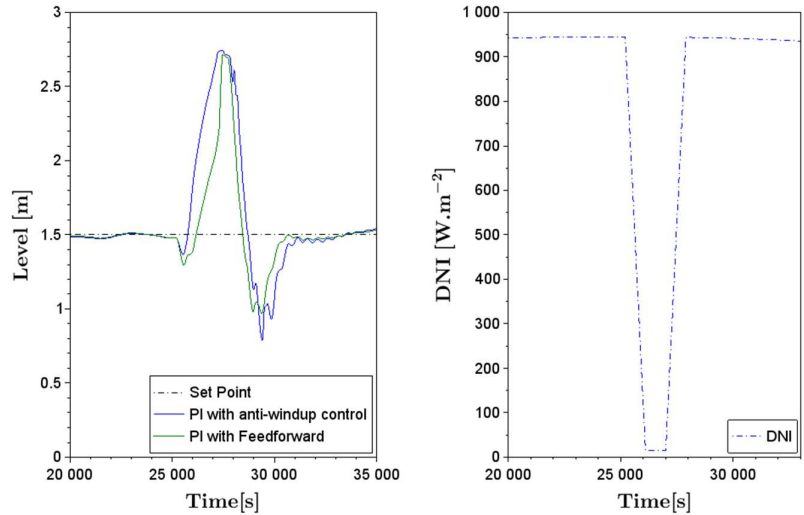

*Figure 4: Evolution of the fluid level (left) during a 15 minutes cloud passage transient (right) – numerical results obtained from the CSP multi-domain simulator.*

We are currently working on validating the simulator against experimental data for a wide range of weather conditions and operation strategies. In a near future, the simulator will be used to elaborate a coherent set of safety functions for the CSP plant.

## 4.2. Development of an advanced control strategy for a multi-energy production plant

The second application presents the methodology adopted for developing an advanced control strategy for a multi-energy production plant. Multi-energy production plants combining various energy vectors, including renewable ones, are a promising direction for more efficient and more sustainable district heating networks and smart energy systems in general [17]. As an example, the potential of a production plant combining a biomass boiler, a heat pump, and a heat storage in three different European contexts has been studied in previous work [18]. However, the operational control required for such a plant to perform efficiently is complex. The PEGASE platform facilitates the development and validation of an advanced Energy Management System (EMS), by combining the co-simulation and control features.

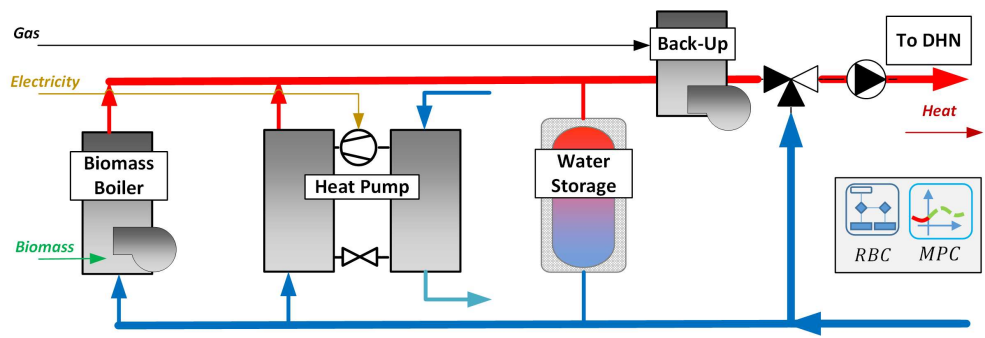

*Figure 5: Overview of the multi-energy DH production plant controlled by the EMS.*

Fig. 5 gives an overview of the plant we consider. To design the EMS, we first developed a numerical simulator of this plant by combining pre-existing models from the Modelica *DistrictHeating* library [19]. The numerical simulator is encapsulated as an FMU. We then programmed an optimization problem using the MILP modeller described in section 3.2. This optimization problem consists of an objective function representing the plant production cost and constraints accounting for, among others, technical limitations of the components and load satisfaction.

We used the PEGASE platform to couple these two components. At each time step, the control problem is formulated and solved using forecasted boundary conditions and system states provided by the numerical simulator. Conversely, the set points periodically updated by the EMS feed the numerical simulator. We typically consider a 15 minutes time step and a rolling horizon of 24 hours, leading to MILP optimization problems comprising around 700 continuous variables, 360 binary variables and 2400 constraints. With these settings, yearly simulations, which involve solving over 30.000 such optimisation problems, last less than 3 hours on a conventional laptop.

This advanced EMS is compared to standard Rule-Based Control (RBC), also implemented in PEGASE using the graphical logic editor. Fig. 6 shows that an advanced EMS naturally benefits from electricity prices variations and thermal storage capacity, while this effect is difficult to capture using RBC. In the analysed context, we demonstrate a 5% reduction in the operational costs when the EMS is used, which could amount to several thousand euros of yearly savings in a real DH network.

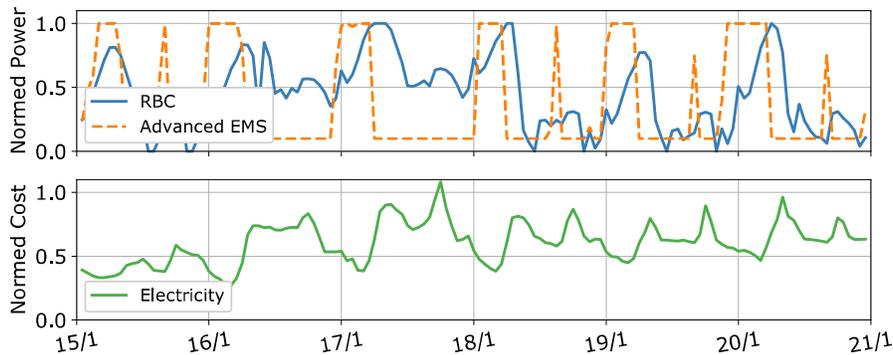

*Figure 6: Advanced EMS and RBC performance to control a multi-energy DH plant.*

### 4.3. Deployment of an optimal controller for a phase change material heat storage

In the frame of the City-Zen FP7 European demonstration project, CCIAG, the DH operating company in the city of Grenoble, France, is currently building an innovative thermal loop in the newly-built Flaubert district [20]. This network is coupled to the main city DH network with a High Pressure / Low Pressure (HP/LP) heat exchanger. As shown in Fig. 7, this system features a phase change material thermal energy storage module and a thermal solar field.

Within this project, we used the PEGASE platform to implement an advanced EMS, whose objective is to optimize the heat power provided by the main heating network with respect to production cost and carbon content. Indeed, typical DH networks feature morning and evening heat consumption peaks, due to end-users behavior. This may result in using costly and polluting peak generators.

Each component (e.g. phase change material heat storage) is regulated by low-level Rule-Based Controls (RBCs), whereas the operating strategy of the system is defined by an advanced EMS. To develop this EMS, we used the same approach as the previous application. Furthermore, we connected the EMS to the real system for operational control. The EMS is installed in parallel to the standard DH operator control system and feeds the Programmable Logic Controller (PLC) with external control set points. By using a dedicated remote control signal, the operator can then choose between the EMS set point and the standard set point. Real time measurements are also sent by the PLC to the operator SCADA through a remote signal. The measurements data feed a database, and the EMS uses those updated values to compute optimized control set points trajectories on a real-time basis (each 15 minutes in this case).

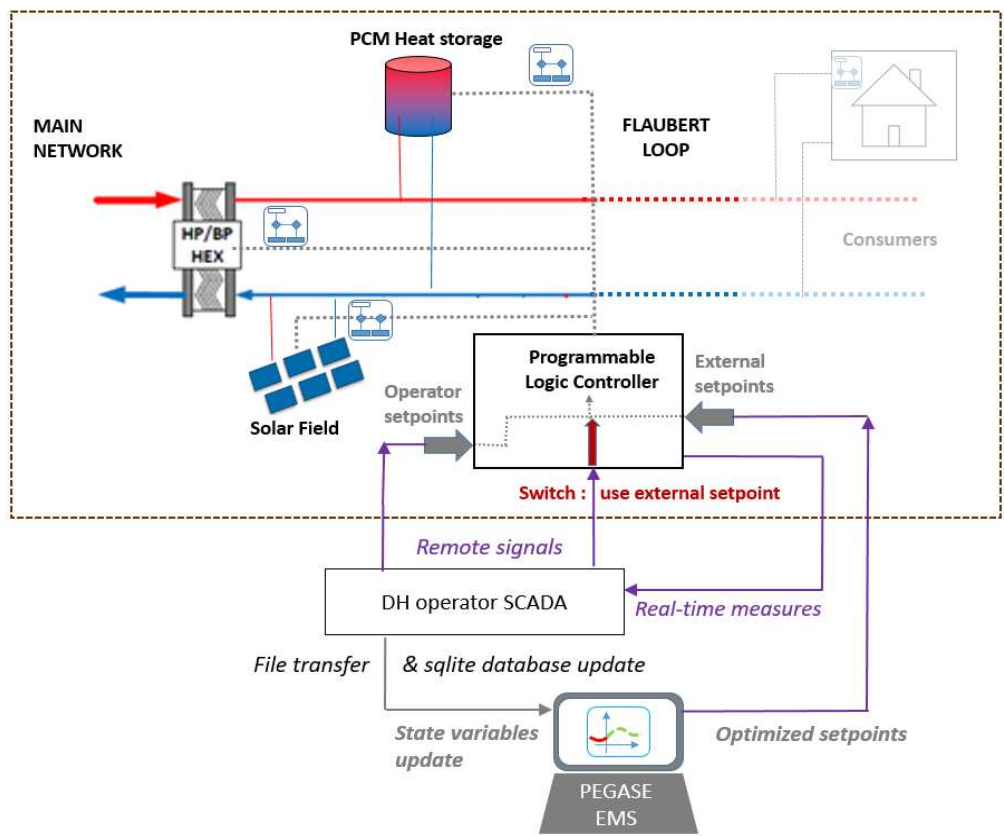

*Figure 7 : The production system of the thermal loop and its data acquisition and control system.*

## 4.4. Techno-economic assessment of a thermal storage in a DH production plant

The co-simulation and control capabilities are also useful for performing long-term techno-economic assessment studies. In particular, MPC control strategies are beneficial when storage is involved in the system, and taking such control into account influences the optimal sizing of components.

In the application presented in this section, we consider a district heating (DH) production system with a biomass boiler, a gas boiler and a thermal storage. This production mix is typical of small-to-medium scale DH network in France, which provide a high share of renewable heat from biomass. However, a biomass boiler is a costly equipment and its size is often optimized to meet a renewable share target. It is also much less flexible than a gas boiler and must operate at least at a given percentage of its nominal power. Including a thermal storage together with optimal control has the potential of increasing the use of the biomass boiler and renewable share, at lower cost.

In this study, we used PEGASE to perform yearly simulations of the DH production system, with an MPC calculation performed every 24h over a 48h rolling horizon. It can be noticed that rolling horizon may not be needed in this kind of study, in which all input data is known in advance and could be integrated into a single optimisation problem. This, however, would result in an optimisation problem too large to be solved in a reasonable time by typical solvers. Instead, a rolling horizon approach allows solving the problem in chunks, taking advantage of the fact that coupling is limited [21]. Indeed, a decision taken at one point mostly influences the following hours of operation, but not the whole year.

To perform this study, a typical heat load curve for a district of 1500 housings was used, based on historical weather data. This load curve is a time series of one year with a time step of one hour. In order to keep the resolution time short, the modelling details are low: only power flows are considered and the thermal storage is modelled as an energy box. Start and stop of the biomass boiler are modelled using a binary variable, with constraints on a minimum duration between two stops of the biomass boiler. The constraint on minimum biomass boiler production level is also taken into account. Details on the parameters are given in Tab. 1. It should be noted that we are considering a thermal storage with a short minimum discharge time (2 hours).

*Table 1: Parameters of the techno-economic study*

| Parameter | Value | Parameter | Value |
|---|---|---|---|
| Annual energy load | 21 217MWh | Natural gas price | 35€/MWh HHV |
| Renewable share target | 60% | Biomass price | 30€/MWh HHV |
| Gas boiler maximum power (load peak power) | 9.8MW | Biomass boiler maximum power | 3.05MW |
| Thermal storage capacity | 2MWh | Biomass boiler minimum power | 40% of Pmax |
| Thermal storage maximum charge/discharge power | 1MW | Minimum duration between 2 stops of biomass boiler | 10 hours |

The MILP problem is solved with a time step of 1 hour 24 hours per 24 hours with a rolling time horizon of 48 hours. A perfect knowledge of the load is assumed. The solution provides the optimal production level for the 2 boilers and the optimal use of thermal storage hour by hour over the year considered. The results confirm that the share of the biomass boiler is increased by 8% over the year using the considered thermal storage (from 12 841 MWh to 14 464 MWh).

## 5. Summary and perspectives

In this paper, we presented a novel research software, called PEGASE, suitable for the design, validation and deployment of advanced control strategies for complex multi-domain energy systems. This platform has been developed based on two requirements identified from the state-of-the art regarding co-simulation and control of complex energy systems : (i) providing co-simulation support accessible to non-expert users, especially regarding the implementation of the sub-systems coupling, and (ii) enabling two types of control: Rule Based Control (RBC) and Model-Predictive Control (MPC). To this end, PEGASE especially features a highly efficient co-simulation engine, together with integrated solutions for defining both RBC and MPC control strategies.

The overall architecture and generic components of the PEGASE platform are presented, with a focus on their distinctive features. In order to illustrate its capabilities, we also present four different applications built on top of this platform. The first application focuses on co-simulation combining two dedicated modelling tools, in the field of Concentrated Solar Power (CSP) plants. The second application presents the generic methodology adopted for designing and validating the designed of advanced Energy Management Systems (EMS). The third application illustrates the deployment of such an EMS on a real-world demonstrator. Finally, the fourth application introduces the use of these methodologies for optimal sizing of an energy system. The variety of examples demonstrates the robustness and genericity of the approach. Further work includes the application of this platform to several domains, as well as the increased usage in the context of techno-economic assessment and system sizing taking into account significant levels of uncertainties in parameters and in predictions.

## Acknowledgments

Work on the PEGASE platform has been conducted in the frame of several projects, including but not limited to [CityZen](FP7 n°608702) and [PENTAGON](H2020 n°731125). Financial support of CCIAG is also acknowledged.

## Nomenclature

**Acronyms**

| | |
|---|---|
| API | Application Programming Interface |
| CCIAG | Compagnie de Chauffage Intercommunale de l'Agglomération Grenobloise |
| CEA | French Alternative Energies and Atomic Energy Commission |
| CSP | Concentrated Solar Power |
| DH | District Heating |
| DSG | Direct Steam Generation (CSP) |
| EMS | Energy Management System |
| FBSF | Functional Bloc Simulation Framework |
| FMI | Functional Mock-up Interface |
| FMU | Functional Mock-up Unit |
| HP/LP | High Pressure / Low Pressure |
| MILP | Mixed Integer Linear Programming |
| MPC | Model Predictive Control |
| PtoH | Power To Heat |
| PLC | Programmable Logic Controller |
| RBC | Rule Based Control |
| RES | Renewable Energy Sources |
| SCADA | Supervisory Control And Data Acquisition |